\documentclass{amsart}

\usepackage{xcolor}
\newtheorem{theorem}{Theorem}[subsection]
\newtheorem{lemma}[theorem]{Lemma}
\theoremstyle{definition}
\newtheorem{definition}[theorem]{Definition}

\newtheorem{corollary}[theorem]{Corollary}
\newtheorem{proposition}[theorem]{Proposition}
\theoremstyle{remark}
\newtheorem{remark}[theorem]{Remark}

\usepackage{enumerate}
\usepackage{amsmath}
\begin{document}
\def\C{\mathcal C}
\def\R{\mathbb R}
\def\X{\mathbb X}
\def\cA{\mathcal A}
\def\cT{\mathcal T}
\def\A{\mathbb A}
\def\B{\mathcal B}
\def\Z{\mathbb Z}
\def\P{\mathbb P}
\def\I{\mathbb I}
\def\H{\mathbb H}
\def\Y{\mathbb Y}
\def\Z{\mathbb Z}
\def\N{\mathbb N}
\def\cal{\mathcal}
\def\ve{\varepsilon}
\def \dis {\displaystyle}
\def \dint{\displaystyle{\int}}
\def \dintt{\dint_{\!\!\!0}^t}
\def\r{\rightarrow }
\newcommand \bp {{\textbf{\bf Proof. }}}
\newcommand{\cqfd}
{\mbox{}\nolinebreak\hfill\rule{2mm}{2mm}\medbreak\par}
\title[]{
The generalization of   Schr\"oder's theorem(1871):\break ~~The multinomial theorem for formal power series under composition }

\author{Galamo  Monkam}  \noindent \thanks{Galamo  Monkam, Department of Mathematics, Morgan State
University, 1700 E. Cold Spring Lane, Baltimore, 21251, USA.\ \\
Email: galamo.monkam@morgan.edu; mondelfi2001@yahoo.fr}
\maketitle

\noindent \textbf{Abstract:}~
\rm {We consider formal power series $ f(x) = a_1 x + a_2 x^2 + \cdots $ \ $(a_1 \neq 0)$, with 
\break coefficients in a field.\ 
We revisit the classical subject of iteration of formal power series,~ the n-fold composition  $f^{(n)}(x)=f(f(\cdots (f(x)\cdots))=f^{(n-1)}(f(x))=\sum\limits_{k=1}^{\infty}f_k^{(n)}x^k$ for $n=2,3,\ldots$,\ where  $f^{(1)}(x)=f(x)$.\ The study of this was begun, and the coefficients $f_k^{(n)}$ where calculated  assuming $a_1=1$, by Schr\"oder  in 1871.\ 
The major result of this paper, Theorem \ref{th41},~generalizes   Schr\"oder \cite{Schr1871}.~It gives explicit formulas for the coefficients $f_k^{(n)}$ when $a_1 \neq 0$ and it is viewed as an analog to the $n^{th}$ Multinomial Theorem Under Multiplication.~We prove Schr\"oder's Theorem using a new and shorter approach.~The Recursion Lemma, Lemma \ref{t1},
which sharpens Cohen's lemma (Lemma \ref{l22}) \ 
is our key tool in the  proof of Theorem \ref{th41} and it  is one of  the most useful recurrence relations for the composition of formal power series.~Along the way we develop numerical  formulas for ~ $f_k^{(n)}~ {\rm where}~1\leq k\leq 5$. \\



\section{\bf Introduction} 
\par \rm  In 1871, the problem of finding the nth iterative power of a formal power series,~$f^{(n)}(x)$,  was introduced and solved by the famous mathematician Ernst  Schr\"oder  in the case $a_1 = 1$~\cite{Schr1871}.  He left the general case open saying that "the formula when $a_1 = 1$ can be found much more easily than in the general case".~He thus asserted that the case $a_1 \neq 1$ was much harder, and this has not been done till now.
\par A main challenge in Mathematics is to give closed formulas for quantities which are defined recursively.
A formal power series {\textit f(x)}~is an expression of the form 
$f(x)=a_0+a_1x+a_2x^2+\cdots $. The word "formal" indicates that this is actually just the sequence of its coefficients. The composition of $f(x)$ and $g(x)$ is defined by \[(f\circ g)(x)=f(g(x))=a_1\cdot g(x)+a_2\cdot(g(x))^2+\cdots +a_n(g(x))^n+\cdots\]The condition $a_0=0$ insures that the coefficient of $x^n$ in $f(g(x))$ is a finite sum. In this paper we assume $a_0=0$. If $f(x) = a_1x + a_2x^2+\cdots$, our main purpose is to study the recursively defined   n-fold composition  by,\ ~$f^{(1)}(x)=f(x)$ and \[f^{(n)}(x)  = f(f(\cdots f(x))\cdots) = f^{(n)}_1 \cdot x + \cdots + f^{(n)}_k \cdot x^k +\cdots ~ \quad \forall n, ~n~>1.\]
 \par {\rm \quad This paper investigates the coefficient $f_k^{(n)}$  for almost unit ($a_1\neq 0$) formal power series} {\it f}  {\rm and  improves the results obtained by  Schr\"oder \cite{Schr1871} and Cayley \cite{Cayley} by giving a recursive and  non-recursive formulation of the composition of formal power series whenever $a_1\neq 0$.}~The major result of this paper, Theorem \ref{th41}, provides a non-recursive formula for $f^{(n)}_k $, for $k\geq 1$ and 
$n\in \Bbb N$. This major result  may be viewed as an analogue of the well-known Multinomial 
Theorem ~(Theorem \ref{t0}),~which gives the coefficients $a_k^{[n]}$~of the product $\big(f(x)\big)^n$.

\section{\bf  Background on Formal Power Series.}
\subsection{Preliminaries}
\ \\
\begin{definition}~{\rm(~From \cite{gk02}~)}
\ \\ \rm \quad Let $S$ be a ring,~a formal power series in one variable \rm on $S$ is defined to be a mapping
from $ {\Bbb N} $ to $S$,~where $\Bbb N$ represents the natural numbers.~We denote the set of all such mappings by $ {\Bbb X} (S) $, \ or $ \Bbb X $. A formal power series $f$ in $x$
from $ \Bbb N$ to $S$ is usually denoted by \ \\
\centerline{ $ f(x) =  a_1 x + \cdots + a_n x^n + \cdots $, \ where \ 
$ \{ a_j \}_{j=1}^\infty \subset S $.} \medskip
In this thesis, we only consider formal power series 
 $f(x)\in {\Bbb X}$.
\end{definition}
\begin{definition}\ \\
{\rm If}  $f(x)\ {\rm and }~ g(x)$\ {\rm are two formal power series} {\rm such that:}
\[f(x)=a_1x+a_2x^2+\cdots \ {\rm and} \quad  g(x)=g_1x+g_2x^2+\cdots\] 
{\rm then} \ \\
{\rm (1)}\quad  {\rm The composition of} $f(x)$ {\rm and} $g(x)$ {\rm is defined by:} \[(f\circ g)(x)=f(g(x))=a_1\cdot g(x)+a_2\cdot(g(x))^2+\cdots +a_n(g(x))^n+\cdots\]
\rm The condition that $g_{0} =0$ insures that the coefficient of $x^n$ in $(f\circ g)(x)$ is finitely calculable as the coefficient of $x^n$ in $\sum\limits_{i=1}^{n}a_i(g(x))^i$.\ \\
{\rm (2)}\quad {\rm The multiplication of } $f(x)$ {\rm and} $g(x)$ {\rm is defined by:}
\qquad \[ (f \cdot g ) (x) = g(x) \cdot f(x) = \sum_{n=1}^\infty \, c_n \, x^n , \ c_n = \sum_{j=1}^n g_j a_{n-j} , \ n =  1, 2, \cdots \] . 
\end{definition}
\subsection{The Multinomial Theorem for Formal Power Series Under \break Multiplication}\label{s1.1}\ \\
\par \rm We first investigate the coefficients of $ (f (x))^{i} $ if $ f(x) $ is
a formal power series. ~Mathematical induction or the multinomial 
coefficients can be used to initiate the investigation of the coefficients of \ 
$(f(x))^{i} $. 
\begin{definition} \ \\
\rm  Let \ $ f(x) =a_1 x + a_2 x^2 + \dots + a_k x^k + \dots $ \ be a formal
power series. \ \ For every \ $ i \in \Bbb N $, \ we write \\ \medskip
{ $ (f(x))^i =\underbrace{(f(x)\cdot~ f(x)\cdot~ \cdots ~\cdot f(x))}_{i~times}$
\par
$\hskip 0.2in  = a_1^i x^i + ia_1^{i-1}a_2 x^{i+1} +
\dots = a_i^{[i]} x^i +a_{i+1}^{[i]}x^{i+1}+\dots +a_k^{[i]}x^k+\cdots $ }
\noindent  
\ \\ Thus $a_k^{[i]}$ \ is  the coefficient \rm of $x^k ~in $ $ (f(x))^i $.
Note: $a_k^{[i]}=0$ for $k<i$~and ~ \ \\ $ a_k^{[1]} = a_k \quad {\rm for~ all} \ k \in {\Bbb N} $.\\
\end{definition}
\begin{theorem}{\rm The Multinomial Theorem~}\label{t0}
\ \\
\rm
\begin{enumerate}
\item
\rm If f(x) is defined as previously then for any \ $ i ,k  \in \Bbb N $. \[a_k^{[\,i\,]} = \sum_{(r_1, \ldots, r_k)}\frac{i\,!}{(r_1!)(r_2!)\cdots (r_j!) \cdots (r_k!)} \,a_1^{r_1}a_2^{r_2}\cdots a_j^{r_j}\cdots a_k^{r_k}\]
\noindent \rm{where the sum is taken for all possible nonnegative integers \ $ r_1, \cdots, r_k $ \ such that \ $ r_1 +r_2+ \cdots+r_j+\cdots + r_k = i$ \ and \ $ r_1 +
2 r_2+\cdots + jr_j + \cdots + k r_k = k $.
}
\\
\item \rm If $i=2$,~ then for $k\geq 3 $ we have \[a_k^{[2]}=\sum_{i=1}^{k-1}a_ia_{k-i}\]\\
\item $a_k^{[i]}$ is a polynomial in the variable $\{a_1,a_2,\cdots,a_{k-i+1}\}$ \rm ~for ~$k\geq i\geq 2.$
\end{enumerate}
\end{theorem}
\subsection{\bf  History of the Composition of Formal Power Series.}\label{s2.1}
\vskip 0.1in \ \\
\par The composition of formal power series, or functional composition, is a significant feature of formal power series which generates a lot of mathematical
results. Some popular applications of formal series and the composition of them can be seen in composition of analytic functions,~in Riordan Arrays,
in numerical analysis, in differential equations, {\it etc.} (see \cite{hen}, \cite{shw}, and \cite{par}). 

\ \\
\par Schr\"oder proposed that iteration of any function could be studied by iterating a simple conjugate of the function:~ \[f=\phi~ \circ~ g ~\circ ~\phi^{-1} \implies f^{(n)}=\phi ~\circ ~ g^{(n)}~\circ~ \phi^{-1}\].\break ~Scheinberg,~ Schwaiger, and  Reich used the method proposed by Schr\"oder.~Scheinberg  \cite{Schein1970} classified formal power series up to conjugation and used the notion of conjugacy to study the composition of formal fower series. Schwaiger \cite{Schwai} and Reich \cite{Reich} used conjugacy to study roots of formal power series .~Muckenhoupt \cite{Mucken} use direct computation to compute 
the exact expression of $ f_2^{(n)}$ {\rm when }$a_1\neq 1.$
\ \\
\par \rm Cohen (\cite{Cohen2018a},~Lemma 2.4) below is the first for having written the n-fold compositional power  ~$f^{(n)}=\sum\limits_{k=1}^{\infty}f_k^{(n)}x^k$~~of $f(x)=\sum\limits_{m=1}^{\infty}a_mx^m$~~, when $a_1 \neq 1 $,~ so as to give the coefficient of $a_k$ in the polynomial $f_k^{(n)}=f_k^{(n)}(a_1,a_2,\cdots, a_k).$
\subsection{\bf Early Results on the n{\it th} Power of a Series Under Composition.}
\subsubsection{\bf The Lemmas of Muckehoupt (1961) and Cohen (2006) } \ \\

\par We state two lemmas which began the direct computation of $f_k^{(n)}.$
\rm The next chapter~(Chapter 3) generalizes Lemma \ref{l22}, it also provide a new proof of Cohen's Lemma (~\cite{Cohen2018a}, Lemma 2.4~) by providing a complete formulation of his $P_k^{(n)}.$
\ \\
\begin{lemma}{\label{l21}}{Muckenhoupt's Lemma}~(\it \cite{Mucken},Lemma 8~)\ \\
{\rm Given a series }~$f(x)=\sum\limits_{i=1}^{\infty}a_ix^i$,~with $a_1\neq 1$\ then,
\[f^{(n)}(x)=a_1^n \cdot x~+~\frac{a_2(a_1^{2n}-a_1^n)}{a_1^2-a_1} \cdot x^2~+\cdots \cdots\]
\end{lemma}
\begin{lemma}{\label{l22}}{\it Cohen's Lemma (\cite{Cohen2018a}, Lemma 2.4)}\rm \ \\ If $n, k$ are positive integers then there exists a polynomial \ $P^{(n)}_k(x_1, \ldots, x_{k-1})$ \  in k-1 variables, with integer coefficients, such that
\ $P^{(n)}_k(x_1, 0, 0, \ldots,\, 0) = 0$ (i.e., each summand
contains an $x_j$ with $j\neq 1$),\ and such that $P^{(n)}_k$ satisfies the following:
\medskip
\mbox{} \qquad If $f(z) = (a_1z + a_2z^2 + \cdots a_{k }z^{k} + \cdots\ )$ is a formal power series with $a_1\neq 0,$ then
$$f^{(n)}_k  = a_ka_1^{n - 1}\bigg(1 + a_1^{k - 1} +
\cdots\  a_1^{(k-1)(n-1)}\bigg) +  P^{(n)}_k(a_1, a_2, \ldots, a_{k-1})$$
\end{lemma}
\medskip
{\bf Proof:} Please refer to \quad \cite{Cohen2018a} \qquad $\Box$.
\section{\bf The Main Recursive Lemma}
~In this section, we improve Cohen's lemma (Lemma \ref{l22}) by giving the complete formulation of $P_k^{(n)}$ which were unknown.
\subsection{\rm Statement and Proof.}
\begin{lemma}{\rm (~\bf Main Recursive Lemma~)}\label{t1}\ \\
\rm Suppose that 
$  f(x) = a_1 x +\cdots + a_k x^k + \cdots\in \Bbb X$,\quad and \ \\ \quad  $ f^{(n)}(x)=f^{(n)}_1\cdot x+\cdots+f^{(n)}_k\cdot x^k+\cdots$ 
  { then~$f^{(n)}_k$~ can be written~as}\ \\
\[ f^{(n)}_k=a_k~\overbrace{a^{n-1}_1\sum_{i=0}^{n-1}a^{i(k-1)}_1}^{C_{k,n}}~+\overbrace{\sum_{i=0}^{n-2}a^{ki}_{1}\cdot\left[\ \sum_{j=2}^{k-1} f^{(n-i-1)}_ja^{[~j~]}_k~\right]}^{P^{(n)}_k}\]
\end{lemma}
\begin{remark}\ \\
\begin{itemize}
\item $f_k^{(n)}=a_k C_{k,n}+P_k^{(n)}$\ \\
\item $P_k^{(n)}(a_1,0\cdots,0) =0.$ see Cohen (\cite{Cohen2018a}, Lemma 2.4) \  \\
\end{itemize}
\begin{eqnarray*}
\bullet  & C_{k,n}= 
\begin{cases} 
n & ~\rm {if} \ ~ a_1= 1 \\
&\\
0 & ~{\rm if} \ ~ a_1^{k-1}\neq 1~{\rm and}~ a_1^{(k-1)n}=1\\ 
&\\
\frac{a_1^{(k-1)n}-1}{a_1^{k-1}-1}\neq 0,  & ~\rm {if} \ ~a_1^{(k-1)n}\neq  1 \ \ 
\end{cases}
\end{eqnarray*}
\end{remark}\ \\
\noindent {\rm \bf Proof of Lemma \ref{t1}:}\ \\
\ \\ \noindent\rm  By definition ~for $s\geq 2$ we have,\ \\ 
 $f^{(s)}(x)=f^{(s-1)}_1\cdot f(x) +f^{(s-1)}_2\cdot (f(x))^2+\cdots+f_{k-1}^{(s-1)}\cdot
(f(x))^{k-1}+f_k^{(s-1)}\cdot(f(x))^{k}+\dots $ 
\begin{align*}
\noindent \implies &f_k^{(s)}=f^{(s-1)}_1\cdot a_k^{[1]} + f^{(s-1)}_2\cdot a_k^{[2]}+\cdots+ f_{k-1}^{(s-1)}\cdot a_k^{[k-1]} +f_{k}^{(s-1)}\cdot a_k^{[k]}\\ 
  \implies  &f_k^{(s)}-a_1^{k}\cdot f_{k}^{(s-1)} =f^{(s-1)}_1\cdot a_k^{[1]} + f^{(s-1)}_2\cdot a_k^{[2]}+\cdots+ f_{k-1}^{(s-1)}\cdot a_k^{[k-1]} \hskip 0.2in \ ~~ (L_{s-1})
\end{align*}
Taking~a ~linear~ combination ~of~ both~ sides~ of~ these ~equations for $s \in [2,n]$ ~we ~get

\[ \sum_{i=1}^{n-1}L_i(a_1^{k})^{n-i-1}\iff f_k^{(n)}-a_1^{k(n-1)}\cdot f^{(1)}_k=\sum_{i=0}^{n-2}a_1^{ki}\cdot \left[\sum_{j=1}^{k-1}f_j^{(n-1-i)}a_k^{[j]}\right] 
 \]
 \ \\
 so it follows that 
 
\begin{align*}
  f_k^{(n)}=&\ a_1^{k(n-1)}\cdot a_k+\sum_{i=0}^{n-2}a_1^{ki}\cdot\left[f_1^{(n-1-i)}\cdot a_k^{[1]}+\sum_{j=2}^{k-1}f_j^{(n-1-i)}a_k^{[j]}\right] \\
  \\
 =&\ a_1^{n-1}a_k\left[a_1^{(k-1)(n-1)}+\sum_{i=0}^{n-2}a_1^{(k-1)i}\right]+\sum_{i=0}^{n-2}a_1^{ki}\cdot \left[\sum_{j=2}^{k-1}f_j^{(n-1-i)}a_k^{[j]}\right] \\
 \\
 =&\ a_1^{n-1}a_k\left[~\sum_{i=n-1}^{n-1}a_1^{(k-1)i}+\sum_{i=0}^{n-2}a_1^{(k-1)i}\right]+\sum_{i=0}^{n-2}a_1^{ki}\cdot\left[\sum_{j=2}^{k-1}f_j^{(n-1-i)}a_k^{[j]}\right]  \\
 \\
=&\ a_1^{n-1}a_k\sum_{i=0}^{n-1}a_1^{(k-1)i}+\sum_{i=0}^{n-2}a_1^{ki}\cdot\left[\sum_{j=2}^{k-1}f_j^{(n-1-i)}a_k^{[j]}\right] 
\end{align*}
Thus we have proved Lemma \ref{t1}. \qquad $\square$ \ \\
The proof can also be done by induction.
\section{\bf First Computations: {\rm {\rm $f_k^{(n)}$ with k=1,2,3,4,5.}}}
\subsection{Computation of  {\rm $f_k^{(n)}$} with k=1,2,3}
\medskip
\begin{corollary}\label{th1}
\ \\ 
\ \\ {\rm If}~ $ f(x) = a_1 x +a_2x^2 + a_3 x^3+\cdots$,~ and 
~$f^{(n)}(x)=f^{(n)}_1\cdot x+f^{(n)}_2\cdot x^2+f^{(n)}_3\cdot x^3+\cdots$, 
\ then
\begin{align*}
 (a)\ \ f^{(n)}_1=&a^{n}_1\\ 
(b) \ \ f^{(n)}_2 =&a^{n-1}_1a_2(1+a^1_1+\cdots+ a^{n-1}_1)=a^{n-1}_1a_2\sum_{i=0}^{n-1}a^{i}_1\\ 
(c) \ \ f^{(n)}_3 =&a^{n-1}_1a_3\sum_{i=0}^{n-1}a^{2i}_1+2a^{n-1}_{1}a^2_{2}\sum_{i=0}^{n-2}a^{2i}_{1}(1+a_1+\cdots+a^{n-i-2}_1) 
\end{align*}
\end{corollary}
\noindent {\rm \bf Proof~:}\ \\
(b) and (c) are direct consequences of the Lemma \ref{t1}. However, (a),(b), and (c) can be done by induction.
\subsection{\rm Computation of $f_k^{(n)}$ for k=4,5.}
\begin{corollary}{\label{cor1}}{\bf \rm Computation of $ f_4^{(n)}$}\ \\
\rm Suppose that  
\ \\$ f(x) = a_1 x + \cdots +a_4x^4+\cdots$~ and 
~$f^{(n)}(x)=f^{(n)}_1\cdot x+\cdots+f^{(n)}_4\cdot x^4+\cdots$.\\
\noindent If~~$a_1\neq 0$ ~we~ have~:
\ \\
\[f^{(n)}_4=a^{n-1}_1a_4\sum_{i=0}^{n-1}a^{3i}_{1}~+~a^{n-1}_1a_2a_3\left[3a_1\sum_{i=0}^{n-2}a^{3i}_1\cdot\sum_{i_1=0}^{n-i-2}a_1^{2i_1}+2\sum_{i=0}^{n-2}a^{3i}_{1}\cdot\sum_{i_1=0}^{n-i-2}a_1^{i_1}\right]\]
\[\qquad +a^3_2a^{n-2}_1\left[\sum_{i=0}^{n-2}a^{3i}_{1}\cdot\sum_{i_1=0}^{n-i-2}a_1^{i_1}+6a^2_1\sum_{i=0}^{n-3}a^{3i}_1\cdot\sum_{i_1=0}^{n-i-3}a^{2i_1}_1\cdot\sum_{i_2=0}^{n-i-i_1-3}a_1^{i_2}\right]\qquad (1)\] 
\end{corollary}
\noindent {\rm Proof:}\ \\
\rm Lemma \ref{t1}  implies that ~\[f^{(n)}_4=a^{n-1}_1a_4\sum_{i=0}^{n-1}a^{3i}_1+\overbrace{\sum_{i=0}^{n-2}a^{4i}_{1}\cdot\left[\ \sum_{j=2}^{3}f^{(n-i-1)}_ja^{[~j~]}_4~\right]}^{P^{(n)}_4}\]
Using Corollary \ref{cor1} we have 
\[ f^{(n)}_4=a^{n-1}_1a_4\sum_{i=0}^{n-1}a^{3i}_1+a^{[~2~]}_4a_2 a_1^{n-2}\sum_{i=0}^{n-2}a^{3i}_{1}\sum_{i_1=0}^{n-i-2}a_1^{i_1}+a_4^{[3]}a_3a_1^{n-2}\sum_{i=0}^{n-2}a_1^{3i}\sum_{i_1=0}^{n-i-2}a_1^{2i_1}\]
\[\hskip 0.6in +a_4^{[3]}a_3^{[2]}a_2^{[1]}a_1^{n-3}\sum_{i=0}^{n-3}a_1^{3i}\sum_{i_1=0}^{n-i-3}a_1^{2i_1}\sum_{i_2=0}^{n-i-i_1-3}a_1^{i_2}\]
\vskip 0.1in
\noindent By replacing $a_4^{[i]}$ by their respective value for $i\in [1,3]$ the result follows.\\
\newpage
\begin{corollary}{\label{cor2}}{\bf \rm Computation of  $f_5^{(n)}$}\\
 \rm Suppose 
\ \\$ f(x) = a_1 x + \cdots +a_5x^5+\cdots$,~ and 
~$f^{(n)}(x)=f^{(n)}_1\cdot x+\cdots+ f^{(n)}_5\cdot x^5 +\cdots$\,. \\
\rm If~ $a_1\neq 0$, we have:
\begin{align*}
 \quad f^{(n)}_5\ =\, & a^{n-1}_1a_5\sum_{i=0}^{n-1}a^{4i}_{1}+a^{n-1}_1a_2a_4\left[2\sum_{i=0}^{n-2}a^{4i}_1\sum_{i_{1}=0}^{n-i-2}a_1^{i_{1}}+ 4a_1^2\sum_{i=0}^{n-2}a^{4i}_{1}\sum_{i_{1}=0}^{n-i-2}a_1^{3i_{1}}\right]\\ 
&\\
+&\ a^{n-2}_1a_3\Bigg\{3a_1^2a_3\sum_{i=0}^{n-2}a^{4i}_{1}\sum_{i_{1}=0}^{n-i-2}a_1^{2i_{1}}+ a_2^2\Bigg[ 2\sum_{i=0}^{n-2}a^{4i}_{1}\sum_{i_{1}=0}^{n-i-2}a_1^{i_{1}}\\
&\\
+&\ 4a^3_1\sum_{i=0}^{n-3}a^{4i}_1\left( 3a_1\sum_{i_{1}=0}^{n-i-3}a_1^{3i_{1}}\sum_{i_{2}=0}^{n-i-i_1-3}a_1^{2i_{2}}+2\sum_{i_{1}=0}^{n-i-3}a_1^{3i_{1}}\sum_{i_{2}=0}^{n-i-i_1-3}a_1^{i_{2}}\right) \\
&\\
+ &\ 6a_1^2\sum_{i=0}^{n-3}a_1^{4i}\sum_{i_{1}=0}^{n-i-3}a^{2i_{1}}_1\sum_{i_{2}=0}^{n-i-i_1-3}a_1^{i_{2}} + 3a_1\sum_{i=0}^{n-2}a^{4i}_{1}\sum_{i_{1}=0}^{n-i-2}a_1^{2i_{1}}\Bigg]~ \Bigg\}\\
&\\
+&\ a_1^{n-2}a_2^4\Bigg[\sum_{i=0}^{n-3}a_1^{4i}\Bigg( 4a_1^2\sum_{i_{1}=0}^{n-i-3}a^{3i_{1}}_1\sum_{i_{2}=0}^{n-i-i_1-3}a_1^{i_{2}}+ 6a_1\sum_{i_{1}=0}^{n-i-3}a^{2i_{1}}_1\sum_{i_{2}=0}^{n-i-i_1-3}a_1^{i_{2}}\Bigg)\\
\mbox{} & \\
+&\ 24a_1^4\sum_{i=0}^{n-4}a_1^{4i}\sum_{i_{1}=0}^{n-i-4}a_1^{3i_{1}}\sum_{i_{2}=0}^{n-i-i_1-4}a^{2i_{2}}_1\sum_{i_{3}=0}^{n-i-i_1-i_2-4}a_1^{i_{3}}\Bigg]~.
\end{align*}
\end{corollary}
\noindent {\rm Proof:}\ \\
We will provide the proof as an application of the generalization of Schr\"oder's theorem known as the Multinomial theorem for Formal Power Series under Composition, please refer to section \ref{sec43}.~To prove Corollary \ref{cor2}, we can also rely on the proof of Corollary \ref{cor1}.
\section{\rm The general  formula  for $f_k^{(n)}$ for $1\leq k\leq 5$~and~$a_1\neq 0$}
\noindent We summarize Corollary \ref{th1}, Corollary \ref{cor1}, and Corollary \ref{cor2} with the following theorem.

\begin{theorem}\label{th2}{{\rm Computation for $f_k^{(n)}$~ for $1\leq k\leq 5$ }}\ \\
 \rm Suppose\ 
$ f(x) = a_1\cdot x +a_2\cdot x^2+a_3\cdot x^3+a_4\cdot x^4  +a_5\cdot x^5+\cdots $,\, with $a_1\neq 0$, then 
\begin{align*}
f^{(n)}(x)=&\ \underbrace{\Bigg\{a_1^n\Bigg\}}_{f_1^{(n)}}\cdot\, x\\
&\\
+&\ \underbrace{\Bigg\{ a^{n-1}_1a_2\sum_{i=0}^{n-1}a_1^{i}\Bigg\}}_{f_2^{(n)}}\cdot\, x^2\\
&\\
+&\ \underbrace{\Bigg\{ ~ a^{n-1}_1a_3\sum_{i=0}^{n-1}a^{2i}_1+2a^{n-1}_1a^2_{2}\sum_{i=0}^{n-2}a^{2i}_{1}\sum_{i_1=0}^{n-i-2}a_1^{i_1}\Bigg\}}_{f_3^{(n)}}\cdot\, x^3\\
&\\
+&\ \Bigg\{ f_4^{(n)}\Bigg\}\cdot\, x^4\\
&\\
+&\ \Bigg\{f_5^{(n)}\Bigg\}\cdot x^5+\cdots\,, 
\end{align*}
\end{theorem}
\noindent \rm where $f_4^{(n)}$ and $f_5^{(n)}$ were previously given in Corollaries \ref{cor1} and  \ref{cor2}.
\section{\bf Combinatorial Lemma} 
\rm \noindent The following combinatorial Lemma will be used in applying our general results~\ \\($a_1\neq 0$) to get Schr\"oder's results ($a_1=1$).
\begin{lemma}\label{l1}
\ \\ \rm If  $n~\rm {and} ~\alpha $ are two positive integers such that $n> \alpha $ and $i_{m} \in \Bbb N \cup \{0\}$\, for\\ 
$m\in [1,\alpha]$, then 

\begin{eqnarray*}
\binom{n}{\alpha}=
\begin{cases} 
\underbrace{\sum_{i_{1}=0}^{n-\alpha}~\quad \sum_{i_2=0}^{n-i_{1}-\alpha}\cdots \sum_{i_{\alpha}=0}^{n-i_{1}-i_2-\cdots-i_{\alpha-1}-\alpha}1}_{\rm \alpha~summation~symbol}~; & ~\text{if} \ ~ \alpha \geq 2 \\
&\\
\sum\limits_{i_{1}=0}^{n-1}1~  & ~\text{if} \ ~\alpha = 1 \ \ 
\end{cases}\ \\
\end{eqnarray*}
\end{lemma}
\noindent{\bf Remark}
\ \\  \noindent \rm  The proof of the Combinatorial Lemma \ref{l1} has been left to the end of this section. To prove Lemma \ref{l1} we first provide a proof of the following theorem. 
\begin{theorem}\label{th3}
\ \\ \rm If   $n$ and $p$  are two positive integers such that $n \geq p \geq 1$, then 
\begin{align*}
\sum_{p=1}^{n}p(p+1)\cdots (p+\alpha-1)&=\frac{1}{\alpha+1}\prod_{p=n}^{n+\alpha}p~\qquad (1)\\
&\\
&=\frac{n(n+1)\cdots(n+\alpha)}{\alpha+1}\,.
\end{align*} 
\end{theorem}
\noindent {\bf Proof.}
\ \\ \rm We proceed by induction on $n$. \\
For $n=1$, \[\sum_{p=1}^{n}p(p+1)\cdots (p+\alpha-1)=\sum_{p=1}^{1}p(p+1)\cdots (p+\alpha-1)=\alpha !\]
\rm On the other hand, for $n=1$ we have 

\begin{align*}
\frac{1}{\alpha+1}\prod_{p=n}^{n+\alpha}p=&~\frac{1}{\alpha+1}\prod_{p=1}^{1+\alpha}p\\
&\\
=&~\frac{1\cdot2\cdots \alpha(\alpha+1)}{\alpha+1}\\
=&~ \alpha !
\end{align*}
\rm Thus, for $n=1$ we have proved (1).\ \\
Now, we suppose that (1) holds for some $n$,\, $n \geq 1$, and we want to prove the following 
\[\sum_{p=1}^{n+1}p(p+1)\cdots (p+\alpha-1)=\frac{1}{\alpha+1}\prod_{p=n+1}^{n+\alpha+1}p\,.\] 
In fact, 
\begin{align*}
\sum_{p=1}^{n+1}p(p+1)\cdots (p+\alpha-1)=& (n+1)(n+2)\cdots (n+\alpha)+ \sum_{p=1}^{n}p(p+1)\cdots (p+\alpha-1)\\
&\\
=&\ \prod_{p=n+1}^{n+\alpha}p+\frac{1}{\alpha+1}\prod_{p=n}^{n+\alpha}p\\
&\\
=&\ \left(1+\frac{n}{\alpha+1}\right)\prod_{p=n+1}^{n+\alpha}p
\end{align*}
\begin{align*}
=&\ \left(\frac{n+\alpha+1}{\alpha+1}\right)\prod_{p=n+1}^{n+\alpha}p\\
&\\
=&\ \frac{1}{\alpha+1}\prod_{p=n+1}^{n+\alpha+1}p\,. ~\qquad~ \square
\end{align*}
\noindent {\rm \bf{ Proof of Combinatorial Lemma \ref{l1}.}}
\ \\ \rm To prove Lemma \ref{l1}, we only need to prove that for $\alpha\geq 2:$
\[\underbrace{\sum_{i_{1}=0}^{n-\alpha}~\quad \sum_{i_2=0}^{n-i_{1}-\alpha}\cdots \sum_{i_{\alpha}=0}^{n-i_{1}-i_2-\,\cdots\,-i_{\alpha-1}-\alpha}1}_{\rm LHS}=\binom{n}{\alpha}\hskip 1.4in (2)\]
To begin the proof we set:
\[P_1 =~ n-\alpha~{\rm and}~P_j=n-\alpha-\sum_{L=1}^{j-1}i_{L} ~~{\rm when}~j\in \left[2,\alpha\right]; \hskip 0.8in (3) \]
\rm (3) can be translated to 
\begin{equation*}
\rm 
\begin{cases}
P_1 &=~ n-\alpha \\
P_2 &=~ n-i_1-\alpha \\
P_3 &=~ n-i_1-i_2-\alpha \\
\vdots & \\
P_{\alpha}&=~ n-i_1-i_2-\cdots -i_{\alpha-1}-\alpha
\end{cases}
\end{equation*}
\[{\rm so},\, (3)\implies P_j=P_{j-1}-i_{j-1}~{\rm when}~j\in \left[2,\alpha\right]\,. \hskip 1.2in (4)\]
Using the previous relation, we have from~(2):
\begin{align*}
{\rm LHS} &= \sum_{i_{1}=0}^{P_1}~\quad \sum_{i_2=0}^{P_2}\, \cdots \,\sum_{i_{\alpha-1}}^{P_{\alpha-1}}\quad\sum_{i_{\alpha}=0}^{P_{\alpha}}1\\
&\\
& = \sum_{i_{1}=0}^{P_1}~\quad \sum_{i_2=0}^{P_2}\, \cdots \,\sum_{i_{\alpha-2}}^{P_{\alpha-2}}\quad \sum_{i_{\alpha-1}=0}^{P_{\alpha-1}}\left(P_{\alpha}+1 \right) \\
& = \sum_{i_{1}=0}^{P_1}~\quad \sum_{i_2=0}^{P_2}\, \cdots \,\sum_{i_{\alpha-2}}^{P_{\alpha-2}}\quad \sum_{P_{\alpha}=0}^{P_{\alpha-1}}\left(P_{\alpha}+1 \right) \\
&\\
& = \sum_{i_{1}=0}^{P_1}~\quad \sum_{i_2=0}^{P_2}\, \cdots \,\sum_{i_{\alpha-3}}^{P_{\small \alpha-3}}\quad \sum_{i_{\small \alpha-2}=0}^{P_{\alpha-2}}\frac{\left(P_{\alpha-1}+1 \right)(P_{\alpha-1}+2)}{2}~,\quad \rm ~by ~using ~Theorem ~\ref{th3} \\
&\\
& = \frac{1}{2}\cdot \sum_{i_{1}=0}^{P_1}~\quad \sum_{i_2=0}^{P_2}\, \cdots \,\sum_{i_{\alpha-3}}^{P_{\alpha-3}}\quad \sum_{P_{\alpha-1}=0}^{P_{\alpha-2}}(P_{\alpha-1}+1 )(P_{\alpha-1}+2)\\
 = &\ \frac{1}{2}\cdot\sum_{i_{1}=0}^{P_1}~\quad \sum_{i_2=0}^{P_2}\, \cdots \,\sum_{i_{\alpha-3}}^{P_{\alpha-3}}\quad \frac{\left(P_{\alpha-2}+1 \right)(P_{\alpha-2}+2)(P_{\alpha-2}+3)}{3}\quad \rm \quad by~ using ~(1)\\
& \\
 = &\ \frac{1}{2}\cdot \frac{1}{3}\cdot \sum_{i_{1}=0}^{P_1}~\quad \sum_{i_2=0}^{P_2}\, \cdots \,\sum_{i_{\alpha-4}}^{P_{\alpha-4}}\quad \sum_{P_{\alpha-2}=0}^{P_{\alpha-3}}
\quad (P_{\alpha-2}+1)(P_{\alpha-2}+2)(P_{\alpha-2}+3)\,.
\\
& \hskip 1.5in \vdots\\
& \rm By~ following ~the~ same ~pattern, ~we ~have
& \\
{\rm LHS} &= \frac{1}{2}\cdot \frac{1}{3}\cdot\, \cdots\, \cdot \frac{1}{\alpha-1}\cdot~\sum_{i_1=0}^{P_1}~(P_2+1)(P_2+2)\cdots (P_2+\alpha-1)\\
=&\ \frac{1}{(\alpha-1)!}~\cdot~\sum_{P_2=0}^{P_1}~(P_2+1)(P_2+2)\cdots (P_2+\alpha-1)\\
&\\
=&\ \frac{1}{(\alpha-1)!}~\cdot \frac{(P_1+1)(P_1+2)\cdots (P_1+\alpha)}{\alpha}\hskip 1in \rm by~ using ~(1)\\
&\\
=&\ \frac{n(n+1) \cdots (n-\alpha+1)}{\alpha!} \rm ~\hskip 1.3in  by~ using~ relation~(3)\\
&\\
=&\ \binom{n}{\alpha}\,. \hskip 1in \square
\end{align*}
\newpage
\subsection{\rm Obtaining  Schr\"oder's Theorem from Theorem \ref{th2} \ \\ when $1\leq k \leq 5$} 
\begin{theorem}\label{th4}{\rm{Application to Schr\"oder's Theorem}}\ \\
\rm Suppose that 
\ \\$ f(x) = 1\cdot x +a_2\cdot x^2+a_3\cdot x^3+a_4\cdot x^4  +a_5\cdot x^5+\cdots,$~ then 

\begin{align*}
f^{(n)}(x)=&\ \underbrace{\Bigg\{1\Bigg\}}_{f_1^{(n)}}\cdot\, x\\
&\\
+&\ \underbrace{\Bigg\{\binom{n}{1}a_2\Bigg\}}_{f_2^{(n)}}\cdot\, x^2\\\
&\\
+&\ \underbrace{\Bigg\{ ~\binom{n}{2}2a_2^2+\binom{n}{1}a_3~ \Bigg\}}_{f_3^{(n)}}\cdot\, x^3\\
&\\
+&\ \underbrace{\Bigg\{ a_4\binom{n}{1}+(5a_2a_3+a^3_2)\binom{n}{2}+6a_2^3\binom{n}{3}\Bigg\}}_{f_4^{(n)}}\cdot\, x^4\\
&\\
+&\ \Bigg\{f_5^{(n)}\Bigg\}\cdot x^5+ \cdots  \\
\end{align*}
\rm with
\[f_5^{(n)}=a_5\binom{n}{1}+\binom{n}{2}\Bigg(5a_2^2a_3+6a_2a_4+3a_3^2\Bigg )+\binom{n}{3}\bigg(10a_2^4+26a_2^2a_3\bigg)+24\binom{n}{4}a_2^4\,.\]
\end{theorem}
\begin{corollary} \ \\
\rm If $a_1=1$, then Theorem \ref{th2} and Theorem \ref{th4} are equivalent.
\end{corollary}
\noindent \rm {\bf Proof.}\ \\
\rm To show the equivalence of Theorems \ref{th2} and \ref{th4}, we only need to prove that when $a_1=1$, the formulation for $f_k^{(n)}$ in Theorem \ref{th2} is identical to that in Theorem \ref{th4}~ for $1\leq k\leq 5$.\\
\rm Case 1:\, $k=1$.\ \\
\rm From Theorem \ref{th2},\, $f_1^{(n)}=a_1^n$\,
 if $a_1=1$, then it is clear that $f_1^{(n)}=1$, which verifies the Schr\"oder's formula for $f_1^{(n)}.$\ \\
\rm Case 2:\, $k=2$.\ \\
From Theorem \ref{th2},\, the expression for $f_2^{(n)}$ is given by\ \\
 \[f_2^{(n)}=a^{n-1}_1a_2\sum_{i=0}^{n-1}a_1^{i}\,.\]
If $a_1=1$,\, then
\begin{align*}
f_2^{(n)}&= a_2\sum_{i=0}^{n-1}~1\\
&\\
&=a_2\cdot\binom{n}{1} \rm ~~by~ using~ Lemma ~\ref{l1}.
\end{align*} 
\rm Thus,\, if $a_1=1$,\, the expression of\, $f_2^{(n)}$\, is the same as the one in the Schr\"oder's paper, see \cite{Schr1871},~Theorem \ref{th4}.\ \\
\rm Case 3:\, $k=3$.\ \\
From Theorem \ref{th2},~ the expression of $f_3^{(n)}$ is given by
\[f_3^{(n)}~=~ a^{n-1}_1a_3\sum_{i=0}^{n-1}a^{2i}_1+2a^{n-1}_1a^2_{2}\sum_{i=0}^{n-2}a^{2i}_{1}\sum_{i_1=0}^{n-i-2}a_1^{i_1}\,.\]
\rm Now, if $a_1=1$, then
\begin{align*}
f_3^{(n)}~&=~ a_3\sum_{i=0}^{n-1}~1~+~2a^2_{2}\sum_{i=0}^{n-2}\quad\sum_{i_1=0}^{n-i-2}~1\\
&\\
&~=~ a_3\sum_{i=0}^{n-1}~1~+~2a^2_{2}\sum_{i=0}^{n-2}\quad \sum_{i_1=0}^{n-i-2}~1\\
&\\
&~=~a_3\cdot\binom{n}{1}~+~2a_2^2\cdot \binom{n}{2}\rm ~~ by ~using~ Lemma ~\ref{l1}.
\end{align*} 
\rm Thus,\, if $a_1=1$\, the expression of $f_3^{(n)}$\, is the same as the one in the Schr\"oder's paper, see \cite{Schr1871},~Theorem \ref{th4}.\ \\
\rm Case 4:\, $k=4.$\ \\
From Corollary \ref{cor1},\, the expression of $f_4^{(n)}$ is given by
\begin{align*}
f^{(n)}_4&=a^{n-1}_1a_4\sum_{i=0}^{n-1}a^{3i}_{1}+a^{n-1}_1a_2a_3\left[3a_1\sum_{i=0}^{n-2}a^{3i}_1\sum_{i_1=0}^{n-i-2}a_1^{2i_1}+2\sum_{i=0}^{n-2}a^{3i}_{1}\sum_{i_1=0}^{n-i-2}a_1^{i_1}\right]\\
&\\
& +a^3_2a^{n-2}_1\left[\sum_{i=0}^{n-2}a^{3i}_{1}\sum_{i_1=0}^{n-i-2}a_1^{i_1}+6a^2_1\sum_{i=0}^{n-3}a^{3i}_1\sum_{i_1=0}^{n-i-3}a^{2i_1}_1\sum_{i_2=0}^{n-i-i_1-3}a_1^{i_2} \right]\,.
\end{align*}
\rm Now, if $a_1=1$, then  
\begin{align*}
f^{(n)}_4&=a_4\sum_{i=0}^{n-1}~1+a_2a_3\left[3\sum_{i=0}^{n-2}\quad \sum_{i_1=0}^{n-i-2}~1+2\sum_{i=0}^{n-2}\quad \sum_{i_1=0}^{n-i-2}~1\right]\\
&\\
& +a^3_2\left[\sum_{i=0}^{n-2}\quad \sum_{i_1=0}^{n-i-2}~1+6\sum_{i=0}^{n-3}\quad \sum_{i_1=0}^{n-i-3}\quad\sum_{i_2=0}^{n-i-i_1-3}~1\right]\\
&\\
&=a_4 \cdot \binom{n}{1}+(5a_2a_3+a^3_2)\binom{n}{2}+6a_2^3 \cdot \binom{n}{3}~~ \rm by ~using~ Lemma~ \ref{l1}.
\end{align*}
\rm Thus,\, if $a_1=1$, the expression for $f_4^{(n)}$ is the same as the one in Schr\"oder's paper, see \cite{Schr1871}, Theorem \ref{th4}.\ \\
\rm Case 5:\, $k=5.$\ \\
\rm If we replace $a_1=1$ in the expression of $f_5^{(n)}$ in Corollary \ref{cor2}, we have the following
\begin{align*}
f^{(n)}_5&=a_5\sum_{i=0}^{n-1}1 + 6a_2a_4\sum_{i=0}^{n-2}\quad \sum_{i_{1}=0}^{n-i-2}1\\
&\\ 
&+a_3\Bigg\{ a_2^2\Bigg[ 5\sum_{i=0}^{n-2}\quad \sum_{i_{1}=0}^{n-i-2}1+26\sum_{i=0}^{n-3}\quad \sum_{i_{1}=0}^{n-i-3}\quad \sum_{i_{2}=0}^{n-i-i_1-3}1\Bigg] +3a_3\sum_{i=0}^{n-2}\quad \sum_{i_{1}=0}^{n-i-2}1\Bigg\}\\
&\\
&+10a_2^4\sum_{i=0}^{n-3}\quad \sum_{i_{1}=0}^{n-i-3}\quad \sum_{i_{2}=0}^{n-i-i_1-3}~1
+24a_2^4\sum_{i=0}^{n-4}~~\sum_{i_{1}=0}^{n-i-4}~~\sum_{i_{2}=0}^{n-i-i_1-4}~~\sum_{i_{3}=0}^{n-i-i_1-i_2-4}1\\
&\\
&=a_5\cdot\binom{n}{1}+\binom{n}{2}\Bigg(5a_2^2a_3+6a_2a_4+3a_3^2\Bigg )+\binom{n}{3}\bigg(10a_2^4+26a_2^2a_3\bigg)+24a_2^4\cdot \binom{n}{4}\,.
\end{align*}
\rm \noindent Thus,\, if $a_1=1$, the expression for $f_5^{(n)}$ above, is the same as the one in Schr\"oder's paper, see \cite{Schr1871}, Theorem \ref{th4}\,. \hskip 0.3in $\square$ \ \\
\section{\bf  Generalization of Schr\"oder's Theorem(1871): The Multinomial Theorem  Under Composition}{}\label{ch5}
\rm The proof of the Main Theorem is at the end of this chapter in section \ref{sec43}. We want the readers to first understand it  through some examples and applications.
 \subsection{\bf The Main Theorem}
\begin{theorem}\label{th41}\ \\
\rm Suppose\  $ f(x) = a_1 x+\cdots + a_k x^k + \cdots \in \Bbb X$,~~and \ \\ $f^{(n)}(x)=f^{(n)}_1\cdot x+\cdots+f^{(n)}_k\cdot x^k+\cdots$, \qquad
~then 
\ \\
\begin{eqnarray*}
 & f_k^{(n)} = 
\begin{cases} 
a_1^n& \text{if} ~\ k = 1 \\
&\\
  a_2a_1^{n-1}\sum\limits_{i=0}^{n-1}a_1^{i}& \text{if}~ \ k = 2\\
&\\
  a_ka_1^{n-1}\sum\limits_{i=0}^{n-1}a_1^{(k-1)i}+\overbrace{\ \sum\limits_{\alpha= 2}^{k-1}~~A_{\alpha,k}\ }^{P_k^{(n)}}\ & \text{if} ~\ k \geq 3 \ \ 
\end{cases}\ \\
\end{eqnarray*}
where $A_{\alpha,k}$ is a polynomial in the variables $a_1,a_2,\cdots,a_{k-\alpha+1}$ which,~ in its compact form, is given by:
\begin{align*}
& \sum \left(~ \left[a_1^{n-\alpha} a_1^{(j_0-1)i_0}\cdots a_1^{(j_m~-1)i_m}\cdots a_1^{(j_{\alpha-1}~-1)i_{\alpha-1}}\right]\cdot\left[a_{j_{\alpha-1}}a_{j_{\alpha-2}}^{[_{j_{\alpha-1}}]}\cdots a_{j_{m-1}}^{[{j_m]}}\cdots a_{j_1}^{[j_2]}a_{j_0}^{[j_1]} \right]~\right), 
\end{align*}
\rm where the sum is taken over all possible sequences $(i_0,i_1,\cdots i_{\alpha-1},j_0,j_1,\cdots,j_{\alpha-1})$\ \\ \rm such that
\begin{align*}
& 0\leq i_0\leq n-\alpha~,\qquad 2 \leq \alpha \leq k-1 \\
& 0\leq i_m \leq n-\alpha - i_0 -i_1-\cdots - i_{m-1}~,\qquad (1\leq m\leq \alpha-1)\\
& j_0=k>j_1>j_2>\cdots >j_{\alpha-1}\geq 2
\end{align*}
\end{theorem}
In expanded form  $A_{\alpha,k}$ is given by~(~we will later drop the parentheses for simplicity.)
\[\sum \left[B_{\alpha,k}\sum_{i=0}^{n-\alpha}\left[a_1^{(k-1)i}\cdot\sum_{i_1=0}^{n-i-\alpha}\left[a_1^{(j_1-1)i_1}\quad \cdots~\cdot\sum_{i_{\alpha-1}=0}^{n-i-i_1-\cdots-i_{\alpha-2}-\alpha }~\left[a_1^{(j_{\alpha-1}-1)i_{\alpha-1}}\right]\cdots\right]~ \right]~ \right],\]
\ \\
{\rm with} $B_{\alpha,k}(j_1,\cdots ,j_{\alpha-1})= a_1^{n-\alpha}\left(a_{j_{\alpha-1}}a_{j_{\alpha-}2}^{[_{j_{\alpha-1}}]}\cdots a_{j_1}^{[j_2]}a_{k}^{[j_1]}\right)$, \ \\
{\rm where}\quad 
$\left(~j_1>j_2>\cdots >j_{\alpha-1}~\right)$
and \ $\{j_1,j_2,\cdots,j_{\alpha-1}\}$~consists of all possible $(\alpha-1)$ element subsets of~$\{2,3,\cdots,k-1\}$.
\begin{remark}{\label{rm42}}\ \\
\rm If we drop the parentheses the expanded form of $A_{\alpha,k}$ can be written as:
\[\sum \left[B_{\alpha,k}\sum_{i=0}^{n-\alpha}a_1^{(k-1)i}\sum_{i_1=0}^{n-i-\alpha}a_1^{(j_1-1)i_1}\quad \cdots~\sum_{i_{\alpha-1}=0}^{n-i-i_1-\cdots-i_{\alpha-2}-\alpha }~ a_1^{(j_{\alpha-1}-1)i_{\alpha-1}}~~\right]~,\]

\rm To give  a better understanding to the reader we will  be using the expanded form of $A_{\alpha,k}.$
 ~The values of  $A_{\alpha,k}$ when $\alpha=2,~3,~4$~is given by:
\begin{align*}
A_{2,k}=&\sum \left[B_{2,k}(j_1)\cdot \sum_{i=0}^{n-2}\left[a_1^{(k-1)i}\cdot\sum_{i_1=0}^{n-i-2}\left[a_1^{(j_{1}-1)i_1}\right]~\right]\quad \right],\qquad j_1~\in~[2,k-1]\\
&\\
=&\sum \left[B_{2,k}(j_1)\cdot \sum_{i=0}^{n-2}\left[a_1^{(k-1)i}\cdot\sum_{i_1=0}^{n-i-2}a_1^{(j_{1}-1)i_1}\right]~\quad \right],\qquad j_1~\in~[2,k-1]\\
&\\
=&\sum \left[a_1^{n-2}a_{j_1}a_k^{[j_1]} \sum_{i=0}^{n-2}a_1^{(k-1)i}\sum_{i_1=0}^{n-i-2}a_1^{(j_{1}-1)i_1}\right],\qquad j_1~\in~[2,k-1]\\
&\\
A_{3,k}=& \sum \left[B_{3,k}(j_1,j_2)\cdot \sum_{i=0}^{n-3}\left[a_1^{(k-1)i}\cdot\sum_{i_1=0}^{n-i-3}\left[a_1^{(j_{1}-1)i_1}\cdot\sum_{i_2=0}^{n-i-{i_1}-3}\left[a_1^{(j_{2}-1)i_2}\right]~\right]~\right]\quad \right] \\
&\\
=& \sum \left[a_1^{n-3} a_{j_2}a_{j_{1}}^{[j_2]}a_k^{[j_{1}]}\sum_{i=0}^{n-3}a_1^{(k-1)i}\sum_{i_1=0}^{n-i-3}a_1^{(j_{1}-1)i_1}\sum_{i_2=0}^{n-i-{i_1}-3}a_1^{(j_{2}-1)i_2}\right] \\
\end{align*}
{\rm where}~ $j_1>j_2$ and \
 $\{j_1,j_2\}$~consists of all possible 2 elements subsets\break of~$\{2,3,\cdots,k-1\}.$\\
\begin{align*}
& A_{4,k}= \sum \left[B_{4,k}(j_1,j_2,j_3)\cdot \sum_{i=0}^{n-4}a_1^{(k-1)i}\sum_{i_1=0}^{n-i-4}a_1^{(j_{1}-1)i_1}\sum_{i_2=0}^{n-i-{i_1}-4}a_1^{(j_{2}-1)i_2}\sum_{i_3=0}^{n-i-i_1-i_2-4}a_1^{(j_3-1)i_3}\right] \\
&\\
=& \sum \left[a_1^{n-4} a_{j_3}a_{j_{2}}^{[j_3]}a_{j_{1}}^{[j_2]}a_k^{[j_{1}]}\sum_{i=0}^{n-4}a_1^{(k-1)i}\sum_{i_1=0}^{n-i-4}a_1^{(j_{1}-1)i_1}\sum_{i_2=0}^{n-i-{i_1}-4}a_1^{(j_{2}-1)i_2}\sum_{i_3=0}^{n-i-i_1-i_2-4}a_1^{(j_3-1)i_3}\right] \\
\end{align*}
{\rm where}\ $j_1>j_2>j_3$ and 
~ $\{j_1,j_2,j_3\}$~consists of all possible $(4-1)$ elements subsets of~$\{2,3,\cdots,k-1\}.$\\
\end{remark}

\section{\bf Applications of the Main Theorem }

\subsection{\bf Computation of $f_k^{(n)}$ for k=3,4}\ \\

\rm Previously we had
\begin{align*}
(1)\qquad f^{(n)}_3 =&a^{n-1}_1a_3\sum_{i=0}^{n-1}a^{2i}_1+2a^{n-1}_{1}a^2_{2}\sum_{i=0}^{n-2}a^{2i}_{1}(1+a_1+\cdots+a^{n-i-2}_1)\\
\ &\\
(2)\qquad f^{(n)}_4=&a^{n-1}_1a_4\sum_{i=0}^{n-1}a^{3i}_{1}+a^{n-1}_1a_2a_3\left[3a_1\sum_{i=0}^{n-2}a^{3i}_1\sum_{i_1=0}^{n-i-2}a_1^{2i_1}+2\sum_{i=0}^{n-2}a^{3i}_{1}\sum_{i_1=0}^{n-i-2}a_1^{i_1}\right]\\
\ & \\
& +~a^3_2a^{n-2}_1\left[\sum_{i=0}^{n-2}a^{3i}_{1}\sum_{i_1=0}^{n-i-2}a_1^{i_1}+6a^2_1\sum_{i=0}^{n-3}a^{3i}_1\sum_{i_1=0}^{n-i-3}a^{2i_1}_1\sum_{i_2=0}^{n-i-i_1-3}a_1^{i_2}\right]
\end{align*}
\noindent{\bf Notation}\ \\
\rm Our main goal here is to recover $f_3^{(n)},~f_4^{(n)}$~using our non-recursive theorem.
\noindent {\rm Proof:}\ \\
\rm (1)\quad From Theorem \ref{th41}  we have 
\[f_k^{(n)}= a_ka_1^{n-1}\sum_{i=0}^{n-1}a_1^{(k-1)i}+\sum_{\alpha= 2}^{k-1}A_{\alpha,k}\]
It follows that for k=3 we have
\[f_3^{(n)}= a_3a_1^{n-1}\sum_{i=0}^{n-1}a_1^{2i}+\sum_{\alpha= 2}^{2}A_{\alpha,3}=a_3a_1^{n-1}\sum_{i=0}^{n-1}a_1^{2i}+A_{2,3}\]

\[A_{2,k}=\sum \left[a_1^{n-2}a_{j_1}a_k^{[j_1]} \sum_{i=0}^{n-2}a_1^{(k-1)i}\sum_{i_1=0}^{n-i-2}a_1^{(j_{1}-1)i_1}\right],\qquad j_1~\in~[2,k-1]\]
\[\implies A_{2,3}=\sum \left[a_1^{n-2}a_{j_1}a_3^{[j_1]} \sum_{i=0}^{n-2}a_1^{2i}\sum_{i_1=0}^{n-i-2}a_1^{(j_{1}-1)i_1}\right],\qquad j_1=2\]

\[\implies A_{2,3}=\left[a_1^{n-2}a_{2}a_3^{[2]} \sum_{i=0}^{n-2}a_1^{2i}\sum_{i_1=0}^{n-i-2}a_1^{i_1}\right],\qquad~where~ a_3^{[2]}=2a_1a_2\]


\[\implies f_3^{(n)}= a_3a_1^{n-1}\sum_{i=0}^{n-1}a_1^{2i}+ 2a_1^{n-1}a_2^2\sum_{i=0}^{n-2}a_1^{2i}\sum_{i_1=0}^{n-i-2}a_1^{i_1}\ \qquad \square\]
\ \\
(2)\quad from  Theorem \ref{th41} we have 
\[f_k^{(n)}= a_ka_1^{n-1}\sum_{i=0}^{n-1}a_1^{(k-1)i}+\sum_{\alpha= 2}^{k-1}A_{\alpha,k}\]
for k =4 we have 
\[(*)\qquad f_4^{(n)}= a_4a_1^{n-1}\sum_{i=0}^{n-1}a_1^{3i}+\sum_{\alpha= 2}^{3}A_{\alpha,4}=a_4a_1^{n-1}\sum_{i=0}^{n-1}a_1^{3i}+A_{2,4}+A_{3,4}\]
from Theorem \ref{th41} we have the expression of $A_{2,k}$ 

\[A_{2,k}=\sum \left[a_1^{n-2}a_{j_1}a_k^{[j_1]} \sum_{i=0}^{n-2}a_1^{(k-1)i}\sum_{i_1=0}^{n-i-2}a_1^{(j_{1}-1)i_1}\right],\qquad j_1~\in~[2,k-1]\]
\rm If $k=4$ we have:
\[A_{2,4}=\sum \left[a_1^{n-2}a_{j_1}a_4^{[j_1]} \sum_{i=0}^{n-2}a_1^{3i}\sum_{i_1=0}^{n-i-2}a_1^{(j_{1}-1)i_1}\right],\qquad j_1~\in~[2,3]\]

\[A_{2,4}= \left[a_1^{n-2}a_{2}a_4^{[2]} \sum_{i=0}^{n-2}a_1^{3i}\sum_{i_1=0}^{n-i-2}a_1^{i_1}\right]+\left[a_1^{n-2}a_{3}a_4^{[3]} \sum_{i=0}^{n-2}a_1^{3i}\sum_{i_1=0}^{n-i-2}a_1^{2i_1}\right]\]

\[ where\qquad a_4^{[2]}=2a_1a_3+a_2^2~;~a_4^{[3]}=3a_1^2a_2\]
from Theorem \ref{th41} we have the expression of $A_{3,k}$ 
\[A_{3,k} =\sum \left[a_1^{n-3} a_{j_2}a_{j_{1}}^{[j_2]}a_k^{[j_{1}]}\sum_{i=0}^{n-3}a_1^{(k-1)i}\sum_{i_1=0}^{n-i-3}a_1^{(j_{1}-1)i_1}\sum_{i_2=0}^{n-i-{i_1}-3}a_1^{(j_{2}-1)i_2}\right] \]
\rm If k=4 we have:
\[A_{3,4} =\sum \left[a_1^{n-3} a_{j_2}a_{j_{1}}^{[j_2]}a_4^{[j_{1}]}\sum_{i=0}^{n-3}a_1^{3i}\sum_{i_1=0}^{n-i-3}a_1^{(j_{1}-1)i_1}\sum_{i_2=0}^{n-i-{i_1}-3}a_1^{(j_{2}-1)i_2}\right], \ \{j_1,j_2\}\subset \{2,3\}\]

 $j_2<j_1$ $\implies (j_1,j_2)=(3,2)$\ \\
\[\implies A_{3,4} =\left[a_1^{n-3} a_{2}a_{3}^{[2]}a_4^{[3]}\sum_{i=0}^{n-3}a_1^{3i}\sum_{i_1=0}^{n-i-3}a_1^{2i_1}\sum_{i_2=0}^{n-i-{i_1}-3}a_1^{i_2}\right]\]

\[\implies A_{3,4}= 6a_1^{n}a_{2}^3\sum_{i=0}^{n-3}a_1^{3i}\sum_{i_1=0}^{n-i-3}a_1^{2i_1}\sum_{i_2=0}^{n-i-{i_1}-3}a_1^{i_2} \] \ \\ \rm By replacing $A_{2,4}$~and~$A_{3,4}$~in (*) the result follows \qquad $\square$
\newpage
\subsection{\bf Computation of $f_k^{(n)}$ for k=5.}\ \\
\ \\
\noindent \rm Our main goal here is to get the expression of $f_5^{(n)}$ previously given in \break Corollary \ref{cor2}~ using our multinomial theorem for composition of formal power series  (Theorem \ref{th41}). We should mention that the expression of $f_5^{(n)}$ given by \break Corollary \ref{cor2} was computed using the recursive  formulation of $f_k^{(n)}$. \ \\
\ \\
\rm From Theorem \ref{th41} we have:
\[f_k^{(n)}= a_ka_1^{n-1}\sum_{i=0}^{n-1}a_1^{(k-1)i}+\sum_{\alpha= 2}^{k-1}A_{\alpha,k}\]
\rm If $k=5$ we have 
\[f_5^{(n)}= a_5a_1^{n-1}\sum_{i=0}^{n-1}a_1^{4i}+A_{2,5}+A_{3,5}+A_{4,5}\ \qquad \ (**)\]
\rm Now let compute $A_{2,5}$.~From Theorem \ref{th41} we have:

\[A_{2,k}=\sum \left[a_1^{n-2}a_{j_1}a_k^{[j_1]} \sum_{i=0}^{n-2}a_1^{(k-1)i}\sum_{i_1=0}^{n-i-2}a_1^{(j_{1}-1)i_1}\right],\qquad j_1~\in~[2,k-1]\]
\ \\
\rm If $k=5$ then $j_1\in \{2,3,4\}$
\  \\
\begin{align*}
A_{2,5}=&\underbrace{\left[a_1^{n-2}a_{2}a_5^{[2]} \sum_{i=0}^{n-2}a_1^{4i}\sum_{i_1=0}^{n-i-2}a_1^{i_1}\right]}_{j_1=2}+\\
&\\
+&\underbrace{\left[ a_1^{n-2}a_{3}a_5^{[3]} \sum_{i=0}^{n-2}a_1^{4i}\sum_{i_1=0}^{n-i-2}a_1^{2i_1}\right]}_{j_1=3}+\underbrace{\left[a_1^{n-2}a_{4}a_5^{[4]} \sum_{i=0}^{n-2}a_1^{4i}\sum_{i_1=0}^{n-i-2}a_1^{j_{3i_1}}\right]}_{j_1=4} 
\end{align*}
\ \\
\rm Now let compute $A_{3,5}$.\ \\~From Theorem \ref{th41} we have the expression of $A_{3,k}$ 
\[A_{3,k} =\sum \left[a_1^{n-3} a_{j_2}a_{j_{1}}^{[j_2]}a_k^{[j_{1}]}\sum_{i=0}^{n-3}a_1^{(k-1)i}\sum_{i_1=0}^{n-i-3}a_1^{(j_{1}-1)i_1}\sum_{i_2=0}^{n-i-{i_1}-3}a_1^{(j_{2}-1)i_2}\right] \]
\rm If k=5 we have:
\[A_{3,5} =\sum \left[a_1^{n-3} a_{j_2}a_{j_{1}}^{[j_2]}a_5^{[j_{1}]}\sum_{i=0}^{n-3}a_1^{4i}\sum_{i_1=0}^{n-i-3}a_1^{(j_{1}-1)i_1}\sum_{i_2=0}^{n-i-{i_1}-3}a_1^{(j_{2}-1)i_2}\right] \]
$\{j_1,j_2\}\subset \{2,3,4\}\quad and\quad j_1>j_2 \implies (j_1,j_2)\in \{(3,2),(4,2),(4,3)\}$
\begin{align*}
A_{3,5}=&\underbrace{\left[a_1^{n-3} a_{2}a_{3}^{[2]}a_5^{[3]}\sum_{i=0}^{n-3}a_1^{4i}\sum_{i_1=0}^{n-i-3}a_1^{2i_1}\sum_{i_2=0}^{n-i-{i_1}-3}a_1^{i_2}\right]}_{(j_1,j_2)=(3,2)} \\
&\\
+&\underbrace{\left[a_1^{n-3} a_{2}a_{4}^{[2]}a_5^{[4]}\sum_{i=0}^{n-3}a_1^{4i}\sum_{i_1=0}^{n-i-3}a_1^{3i_1}\sum_{i_2=0}^{n-i-{i_1}-3}a_1^{i_2}\right]}_{(j_1,j_2)=(4,2)} \\
&\\
+&\underbrace{\left[a_1^{n-3} a_{3}a_{4}^{[3]}a_5^{[4]}\sum_{i=0}^{n-3}a_1^{4i}\sum_{i_1=0}^{n-i-3}a_1^{3i_1}\sum_{i_2=0}^{n-i-{i_1}-3}a_1^{2i_2}\right]}_{(j_1,j_2)=(4,3)} 
\end{align*}
\rm Now let compute $A_{4,5}$.\ \\~From Theorem \ref{th41} we have the expression of $A_{4,k}$ 
\[A_{4,k}= \sum \left[B_{4,k}(j_1,j_2,j_3)\sum_{i=0}^{n-4}a_1^{(k-1)i}\sum_{i_1=0}^{n-i-4}a_1^{(j_{1}-1)i_1}\sum_{i_2=0}^{n-i-{i_1}-4}a_1^{(j_{2}-1)i_2}\sum_{i_3=0}^{n-i-i_1-i_2-4}a_1^{(j_3-1)i_3}\right]\] 
\ \\
$\{j_1,j_2,j_3\}\subset \{2,3,4\}$~and~$j_1>j_2>j_3\implies (j_1,j_2,j_3)=(4,3,2)$\ \\ 
\ \\
$~{\rm So \quad }B_{4,k}(j_1,j_2,j_3)=a_1^{n-4} a_{j_3}a_{j_{2}}^{[j_3]}a_{j_{1}}^{[j_2]}a_k^{[j_{1}]}=a_1^{n-4} a_{2}a_{3}^{[2]}a_{4}^{[3]}a_5^{[4]}$
\ \\
\[\implies A_{4,5}=  \left[a_1^{n-4} a_{2}a_{3}^{[2]}a_{4}^{[3]}a_5^{[4]}\sum_{i=0}^{n-4}a_1^{4i}\sum_{i_1=0}^{n-i-4}a_1^{3i_1}\sum_{i_2=0}^{n-i-{i_1}-4}a_1^{2i_2}\sum_{i_3=0}^{n-i-i_1-i_2-4}a_1^{i_3}\right]\] 
\ \\ \rm By replacing $A_{2,5}$,~$A_{3,5}$~and~$A_{4,5}$~in (**) the result follows \qquad $\square$
\begin{lemma}\label{l43}\ \\
\rm Let $\alpha$, m be non negative integers such that $\alpha \in [2,k-1]$ and $m\in [1,\alpha-1]$.\ \\ If~ $~k=j_0>j_1>\cdots >j_{m-1}>~j_m>\cdots >j_{\alpha-1}\geq 2~$
and  $\{~j_1,\cdots,j_{m-1},~j_m,\cdots,j_{\alpha-1}~\}$\break  consists of all possible $(\alpha-1)$ element subsets of~$\{2,3,\cdots,k-1\}$ then \ \\
$j_{m-1}-j_{m}\leq k-\alpha$ .
\end{lemma}

\noindent {\bf Proof:}\ \\
Since $~j_1>j_2>\cdots >j_{m-1}>~j_m>\cdots >j_{\alpha-1}\geq 2~$~ we have~ $j_1\geq \alpha.$
\begin{equation*}
j_1 \geq \alpha \implies 
\begin{cases}
j_{m-1} &\geq ~\alpha-m+2\\
j_{m} &\geq ~\alpha-m+1\\
\end{cases}\qquad ;~~
j_1 \leq k-1 \implies 
\begin{cases}
j_{m-1} &\leq ~k-(m-1)\\
j_{m} &\leq ~k-m\\
\end{cases}~\cdot
\end{equation*}
\ \\
Thus it follows that\quad  $j_{m-1}-j_{m}\leq k-\alpha$. \qquad $\square$
\begin{remark}\rm (Remark of Lemma \ref{l43})\ \\
There exists sequences for which the value $j_{m-1}-j_m=k-\alpha$ is achieved.~For instance,~ $k,~ \alpha,~ \alpha-1,~ \cdots ,2$~ is a sequence with $j_0-j_1=k-\alpha~.$
\end{remark}
\subsection{\bf Proof of the Main Theorem }{\label{sec43}}\ \\
\ \\
\rm We only need to prove that $P_k^{(n)}$~defined in Theorem \ref{th41} is the same as the one given in Lemma \ref{t1}.\\
From Lemma \ref{t1}, 
\[P_k^{(n)}=\sum_{i=0}^{n-2}a^{ki}_{1}\left[\ \sum_{j_1=2}^{k-1}f^{(n-i-1)}_{j_1}a^{[j_1]}_k~\right]\]\[ ~{\rm where} ~\qquad {\it f_{j_1}^{(n-i-1)}}=a^{n-i-2}_1a_{j_1}\sum_{i_{1}=0}^{n-i-2}a^{i_{1}(j_{1}-1)}_1+\overbrace{\sum_{i_{1}=0}^{n-i-3}a^{j_{1}i_{1}}_{1}\left[\ \sum_{j_2=2}^{j_1-1}f^{(n-i-i_{1}-2)}_{j_2}a^{[j_2]}_{j_1}~\right]}^{P^{(n-i-1)}_{j_1}}\]
\begin{align*}
\implies P_k^{(n)}=&a_1^{n-2}\sum_{i=0}^{n-2}a_1^{(k-1)i}\sum_{j_1=2}^{k-1}\left[a_{j_1}\sum_{i_1=0}^{n-i-2}a_1^{(j_1-1)i_1}~a_k^{[j_1]}\right]+D_2\\
&\\
=&\sum_{j_1=2}^{k-1}\left[a_1^{n-2}a_{j_1}~a_k^{[j_1]}\sum_{i=0}^{n-2}a_1^{(k-1)i}\sum_{i_1=0}^{n-i-2}a_1^{(j_1-1)i_1}\right]+D_2\\
&\\
=& \underbrace{\overbrace{\sum \left[ B_{2,k}(j_1)\sum_{i=0}^{n-2}a_1^{(k-1)i}\sum_{i_1=0}^{n-i-2}a_1^{(j_1-1)i_1}\right]}}^{A_{2,k}}_{j_1\in [2,k-1]}+D_2\,,\\
\end{align*}
 \[{\rm where}~\ D_2\,=\,\sum_{i=0}^{n-3}a_1^{ki}\sum_{j_1=2}^{k-1}\left[\sum_{i_1=0}^{n-i-3}a_1^{j_1i_1}\sum_{j_2=2}^{j_1-1}\left[f_{j_2}^{(n-i-i_1-2)}a_{j_1}^{[j_2]}\right]~a_k^{[j_1]}\right] ~\cdot \]\\
From Lemma \ref{t1},
\[f_{j_2}^{(n-i-i_1-2)}=a^{n-i-i_1-3}_1a_{j_2}\sum_{i_{2}=0}^{n-i-i_1-3}a^{(j_2-1)i_{2}}_1+\sum_{i_{2}=0}^{n-i-i_1-4}a^{j_2i_{2}}_{1}\left[\ \sum_{j_3=2}^{j_2-1}f^{(n-i-i_{1}-i_{2}-3)}_{j_3}a^{[j_3]}_{j_2}~\right]~\cdot\]
It follows that,
\begin{align*}
D_2=&\ a_1^{n-3}\sum_{i=0}^{n-3}a_1^{(k-1)i}\sum_{j_1=2}^{k-1}\left[\sum_{i_1=0}^{n-i-3}a_1^{(j_1-1)i_1}\sum_{j_2=2}^{j_1-1}\left[a_{j_2}\sum_{i_2=0}^{n-i-1_1-3}a_1^{(j_2-1)i_2}a_{j_1}^{[j_2]}\right]~a_k^{[j_1]}\right] +D_3\\
&\\
=&\ \sum_{j_1=2}^{k-1}\sum_{j_2=2}^{j_1-1}\left[a_1^{n-3}a_{j_2}a_{j_1}^{[j_2]}~a_k^{[j_1]}\sum_{i=0}^{n-3}a_1^{(k-1)i}\sum_{i_1=0}^{n-i-3}a_1^{(j_1-1)i_1}\sum_{i_2=0}^{n-i-1_1-3}a_1^{(j_2-1)i_2}\right] + D_3\\
&\\
=&\ \underbrace{\overbrace{\sum \left[B_{3,k}(j_1,j_2)\sum_{i=0}^{n-3}a_1^{(k-1)i}\sum_{i_1=0}^{n-i-3}a_1^{(j_1-1)i_1}\sum_{i_2=0}^{n-i-1_1-3}a_1^{(j_2-1)i_2}\right]}}^{A_{3,k}}_{ \{j_1,j_2\}~\subset ~\{2,3,\,\ldots\, ,k-1\}~ {\rm and}~ j_1>j_2} +\, D_3\,,\\
\end{align*}
where 
\[ D_3\,=\,\sum_{i=0}^{n-3}a_1^{ki}\sum_{j_1=2}^{k-1}\left[\sum_{i_1=0}^{n-i-3}a_1^{j_1i_1}\sum_{j_2=2}^{j_1-1}\left[\sum_{i_2=0}^{n-i-1_1-3}a_1^{j_2i_2}\sum_{j_3=2}^{j_2-1}\left[f_{j_3}^{(n-i-i_1-i_2-3)}a_{j_2}^{[j_3]}\right]a_{j_1}^{[j_2]}\right]~a_k^{[j_1]}\right]~,\]
\ \\
thus,~$P_k^{(n)}$ can be written as $P_k^{(n)}=A_{2,k}+A_{3,k}+D_3$.\ \\
\noindent By following the same pattern, the stopping point of the conjecture is identified when the index of `$i$' takes on the value $k-2$; this implies that: 
\[P_k^{(n)}=A_{2,k}+A_{3,k}+A_{4,k}+\cdots+A_{\alpha,k}+\cdots +D_{k-2},\, {\rm with}\, D_{k-2}=A_{k-1,k}\,. \quad (**)\].
\rm \hskip 2in So,\,  $P_k^{(n)}=\sum\limits_{\alpha=2}^{k-1}A_{\alpha,k}$\,.\ \\

By setting, \quad $ M_{j_{\alpha-1}}=a_{j_{\alpha-1}}\sum\limits_{i_{\alpha-1}=0}^{n-i-i_1-\,\cdots\,-i_{\alpha-2}-\alpha }a_1^{(j_{\alpha-1}-1)i_{\alpha-1}}\,,~{\rm and}~i_0=i,$
\begin{align*}
&A_{\alpha,k}=a_1^{n-\alpha}\sum_{i=0}^{n-\alpha}a_1^{(k-1)i}\sum_{j_{1}=2}^{k-1}\left[\sum_{i_1=0}^{n-i-\alpha}a_1^{(j_1-1)i_1}\sum_{j_2=2}^{j_1-1}\left[\cdots\sum_{j_{\alpha-1}=2}^{j_{\alpha-2}-1}\left[ M_{j_{\alpha-1}}\cdot a_{j_{\alpha-2}}^{[j_{\alpha-1}]}\right]\cdots a_{j_1}^{[j_2]}\right]~a_k^{[j_1]}\right] \\
&\\
&=\sum_{j_{1}=2}^{k-1}\sum_{j_2=2}^{j_1-1}\cdots\sum_{j_{\alpha-1}=2}^{j_{\alpha-2}-1}\left[ a_1^{n-\alpha}\cdot a_{j_{\alpha-2}}^{[j_{\alpha-1}]}\cdots a_{j_1}^{[j_2]}~a_k^{[j_1]}\sum_{i=0}^{n-\alpha}a_1^{(k-1)i}\sum_{i_1=0}^{n-i-\alpha}a_1^{(j_1-1)i_1}\cdots \left[ M_{j_{\alpha-1}}\right]\right] \\
&\\
&=\sum\left[ B_{\alpha,k}(j_1,j_2,\cdots,j_{\alpha-1})\sum_{i=0}^{n-\alpha}a_1^{(k-1)i}\sum_{i_1=0}^{n-i-\alpha}a_1^{(j_1-1)i_1}\cdots \sum\limits_{i_{\alpha-1}=0}^{n-i-i_1-\,\cdots\,-i_{\alpha-2}-\alpha }a_1^{(j_{\alpha-1}-1)i_{\alpha-1}}\right]~\cdot  \\
\end{align*}
{\rm Here}~~ $B_{\alpha,k}(j_1,\ldots ,j_{\alpha-1})= a_1^{n-\alpha}\left(a_{j_{\alpha-1}}a_{j_{\alpha-2}}^{[_{j_{\alpha-1}}]}\cdots a_{j_1}^{[j_2]}a_{k}^{[j_1]}\right)$ \ \\
with $\{j_p\}$  a strictly decreasing sequence for $\forall p ~\in [1,\alpha-1]$~and~ $\{j_1,j_2,\ldots,j_{\alpha-1}\}$\\consisting of all possible $(\alpha-1)$--element subsets of~$\{2,3,\ldots,k-1\}$.~
Finally, $A_{\alpha,k}$~ involves only the variables $a_1, a_2,\ldots,a_{k-\alpha+1}$ due to (Lemma \ref{l43} and Theorem [\ref{t0}, 3]).
\begin{proposition}\ \\
{\rm The set of all the $B_{\alpha, k}$ in Theorem \ref{th41} are together the sum of ~~$2^{k - 2} - 1$  non-constant summands of the form \ \ \ $a_1^{n_1}a_2^{n_2} \cdots a_{k - 1}^{n_{k - 1}}\ \ (n_i\geq 0)$.  }
\end{proposition}
{\bf Proof.}\,  A set of $k - 2$ elements has $2^{k - 2} - 1$ non-empty subsets. $\square$
\section{\bf A new Proof of  Schr\"oder's Formula for $f_k^{(n)}$ when $a_1=1$ }
\setlength{\parindent}{10ex}
 \rm Schr\"oder \cite{Schr1871} states: "The case $a_1\neq 1$~is much harder":\ \\
{\it If a function $f(z)$ is expressed as a Mac-Laurin power series, then the formation law of  its iterates can be found,
 much more easily than in the general case
\noindent if the series has the particular property that it starts with the beginning term "z". In the following, we will concern ourselves with such a function, whose general case is:\\
$f(z)=z+a_2z^2+a_3z^3+\cdots$ \hfill (34)
\rightline{Schr\"oder(1871)}}
\subsection{\bf Statement and Proof of Schr\"oder's Theorem}\ \\

\rm To provide a new proof of Schr\"oder's theorem, we use the following combinatorial lemma which is equivalent to Lemma \ref{l1}:
\ \\ \rm If $n~ {\rm and} ~\alpha $ are two positive integers such that $n> \alpha $ and $i_{m} \in \Bbb N \cup \{0\}$~for\break  $m\in [0,\alpha-1],~i=i_0$, then  
\ \\
\begin{eqnarray*}
\binom{n}{\alpha}=
 & 
\begin{cases} 
\underbrace{\sum_{i=0}^{n-\alpha}~\quad \sum_{i_1=0}^{n-i-\alpha}\cdots \sum_{i_{\alpha-1}=0}^{n-i-i_1-\,\cdots\,-i_{\alpha-2}-\alpha}1}_{\rm \alpha~summation~symbol}~; & ~\text{if} \ ~ \alpha \geq 2 \\
&\\
\sum\limits_{i=0}^{n-1}1~  & ~\text{if} \ ~\alpha = 1 \ \ 
\end{cases}\ \\
\end{eqnarray*}

\begin{theorem}{\rm{(Schr\"oder's Theorem)}}\label{th51}
\ \\ \rm When $a_1=1$, then
\begin{align*}
{\rm if}~ k=1,\quad  {\rm then}~ \quad &f_1^{(n)}=1; \\
{\rm if}~ k=2, \quad {\rm then}~ \quad &f_2^{(n)}=a_2\binom{n}{1}.
\end{align*}
\rm If $k\geq 3$, then
\ \\
 \[\quad f_k^{(n)}= a_k\cdot\binom{n}{1}+\sum_{\alpha= 2}^{k-1}\left[\binom{n}{\alpha}\cdot \sum a_k^{[j_1]} ~a_{j_1}^{[j_2]}~a_{j_2}^{[j_3]}\cdots~  a_{j_{\alpha-2}}^{[j_{\alpha-1}]} a_{j_{\alpha-1}}^{[\,1\,]}\right]~,\]
\rm where the second sum is taken over those systems of integers\, $j_1,j_2,\cdots, j_{\alpha-1}$,\, which consist of all possible ($\alpha-1$)--element subsets of $\{k-1,k-2,\cdots, 3,2\}.$
\end{theorem}
\ \\
\noindent {\rm \bf Proof.}\, {\rm (of Schr\"oder's Theorem)}\ \\
\rm We will use the Combinatorial Lemma \ref{l1} to transform Theorem \ref{th41} into Schr\"oder's Theorem.

\rm  From Theorem \ref{th41}\,:
\begin{align*}
{\rm For}\quad k=1, \quad {\rm then}\quad & f_1^{(n)}=~a_1^n.\\
{\rm For}\quad k=2, \quad {\rm then}\quad & f_2^{(n)}=~a_2a_1^{n-1}\sum\limits_{i=0}^{n-1}a_1^{i}.
\end{align*}
\rm  If $a_1=1$, then
\begin{align*}
{\rm For}\quad k=1, \quad {\rm then}\quad & f_1^{(n)}=~a_1^n=1.\\
{\rm For}\quad k=2, \quad {\rm then}\quad & f_2^{(n)}=~a_2a_1^{n-1}\sum\limits_{i=0}^{n-1}a_1^{i}=a_2\sum\limits_{i=0}^{n-1}1=a_2\binom{n}{1}. 
\end{align*}

\rm \noindent Therefore, when $k \in [1,2]$, the result is transparent. 
\qquad $\square$\ \\
\ \\
\rm Now, from Theorem \ref{th41}, when $k\geq 3:$
\[ f_k^{(n)}= a_ka_1^{n-1}\sum_{i=0}^{n-1}a_1^{(k-1)i}+\sum_{\alpha= 2}^{k-1}A_{\alpha,k}\,,\]

\[~{\rm where}\quad A_{\alpha,k}=\sum \left[B_{\alpha,k}\sum_{i=0}^{n-\alpha}a_1^{(k-1)i}\sum_{i_1=0}^{n-i-\alpha}a_1^{(j_1-1)i_1}\cdots \sum_{i_{\alpha-1}=0}^{n-i-i_1-\,\cdots\,-i_{\alpha-2}-\alpha }a_1^{(j_{\alpha-1}-1)i_{\alpha-1}} \right]\]
\ \\
{\rm with}\quad $B_{\alpha,k}=B_{\alpha,k}(j_1,\cdots ,j_{\alpha-1})= a_1^{n-\alpha}\left(a_{j_{\alpha-1}}a_{j_{\alpha-}2}^{[{j_{\alpha-1}}]}\cdots a_{j_1}^{[j_2]}a_{k}^{[j_1]}\right)$.\ \\

\rm If $a_1=1$, then 
\[f_k^{(n)}=a_k\cdot\binom{n}{1}+\sum_{\alpha= 2}^{k-1}A_{\alpha,k}\,,\]
{\rm where}
\begin{align*}
 A_{\alpha,k}=&\ \sum \left[\left(a_{j_{\alpha-1}}a_{j_{\alpha-}2}^{[{j_{\alpha-1}}]}\cdots a_{j_1}^{[j_2]}a_{k}^{[j_1]}\right)\cdot \sum_{i=0}^{n-\alpha}~~\sum_{i_1=0}^{n-i-\alpha}\cdots \sum_{i_{\alpha-1}=0}^{n-i-i_1-\,\cdots\,-i_{\alpha-2}-\alpha }~1 \right]\\
 &\\
 =&\ \sum \left[\left(a_{j_{\alpha-1}}a_{j_{\alpha-2}}^{[{j_{\alpha-1}}]}\cdots a_{j_1}^{[j_2]}a_{k}^{[j_1]}\right)\cdot\binom{n}{\alpha} \right] \quad {\rm (see~Lemma~\ref{l1}})\\
 &\\
 =&\ \binom{n}{\alpha}\cdot \sum a_{j_{\alpha-1}}a_{j_{\alpha-}2}^{[{j_{\alpha-1}}]}\cdots a_{j_1}^{[j_2]}a_{k}^{[j_1]},\\
\end{align*}
{\rm where}\, $\{j_p\}$ is a strictly decreasing sequence for $\forall\, p ~\in [1,\alpha-1]$ and \ \\
\ $\{j_1,j_2,\cdots,j_{\alpha-1}\}$ consists of all possible $(\alpha-1)$ elements subsets of~$\{2,3,\cdots,k-1\}$.\ \\
\ \\
\rm It follows that,~for $k\geq 3$, if $a_1=1$,~then:
\ \\
\begin{align*}
 f_k^{(n)}=&\ a_k\cdot\binom{n}{1}+\sum_{\alpha= 2}^{k-1}\left[ \binom{n}{\alpha}\sum a_k^{[j_1]} ~a_{j_1}^{[j_2]}\cdot a_{j_2}^{[j_3]}~ \cdots~a_{j_{\alpha-2}}^{[j_{\alpha-1}]} a_{j_{\alpha-1}}^{[1]}\right]\,. \ \qquad \square \\
\end{align*}
\newpage
\subsection{\bf~ Computation of $f_5^{(n)}$ Using Theorem \ref{th51}}\ \\
\ \\
\rm \noindent From Theorem \ref{th51} we have: \ \\
 \rm if~ $k \geq 2 $, then 
 \[\quad f_k^{(n)}= a_k\cdot\binom{n}{1}+\underbrace{\sum_{\alpha= 2}^{k-1}\left[\overbrace{ \binom{n}{\alpha}\cdot \sum a_k^{[j_1]}\cdot ~a_{j_1}^{[j_2]}\cdot a_{j_2}^{[j_3]}\cdot\, \cdots\, \cdot a_{j_{\alpha-2}}^{[j_{\alpha-1}]}\cdot a_{j_{\alpha-1}}^{[\,1\,]}}^{A_{\alpha,k}}\right]}_{P_{k}^{(n)}}~ {\bf \cdot} \]
\rm For~ $k=5$,
\[f_5^{(n)}=a_5\cdot\binom{n}{1}+\sum_{\alpha= 2}^{4}\left[ \underbrace{\binom{n}{\alpha}\cdot \sum a_5^{[j_1]}\cdot ~a_{j_1}^{[j_2]}\cdot a_{j_2}^{[j_3]}\cdot\, \cdots\, \cdot a_{j_{\alpha-2}}^{[j_{\alpha-1}]}\cdot a_{j_{\alpha-1}}^{[\,1\,]}}_{A_{\alpha,5}}\right]~ {\bf \cdot}\]
If $\alpha=2$, then 
\begin{align*}
 A_{2,5} =&\ \binom{n}{2}\cdot\sum a_5^{[j_1]}\cdot a_{j_1}^{[\,1\,]}\quad {\rm where}~j_1\, \in \{2,3,4\}\\
=&\ \binom{n}{2}\cdot\left(a_5^{[\,2\,]}\cdot a_2^{[\,1\,]}+a_5^{[\,3\,]}\cdot a_3^{[\,1\,]}+a_5^{[\,4\,]}\cdot a_4^{[\,1\,]}\right)\\
=&\ \binom{n}{2}\cdot\left( ~5a_2^2+6a_2a_4+3a_3^2 ~\right)\, {\bf \cdot}
\end{align*}
If $\alpha=3$, then 
\rm ~$\{j_1,j_2\}$ consists of all possible 2--element subsets of $\{2,3,4\}$, where $j_2<j_1$, so we have\, 
$\{j_1,j_2\}\in\{ \{4,3\},\{4,2\},\{3,2\} \},$\, {\rm and}
\begin{align*}
 A_{3,5}=&\ \binom{n}{3}\cdot\sum a_5^{[j_1]}\cdot a_{j_1}^{[j_2]}\cdot a_{j_2}^{[\,1\,]}\\
=&\ \binom{n}{3}\cdot\left(a_5^{[\,4\,]}\cdot a_4^{[\,3\,]}\cdot a_3^{[\,1\,]}+a_5^{[\,4\,]}\cdot a_4^{[\,2\,]}\cdot a_2^{[\,1\,]}+a_5^{[\,3\,]}\cdot a_3^{[\,2\,]}\cdot a_2^{[\,1\,]}\right)\\
=&\ \binom{n}{3}\cdot\left(10a_2^4+26a_2^2a_3 \right)~ {\bf \cdot}
\end{align*}
\rm If $\alpha=4$, then 
~$\{j_1,j_2,j_3\}$~consists of all possible 3--element subsets of $\{2,3,4\}$,\ \\ \rm where $j_3<j_2<j_1$. That is, 
$\{j_1,j_2,j_3\}= \{4,3,2\}$. We have,
\begin{align*}
\binom{n}{\alpha}\cdot\sum a_5^{[j_1]}\cdot ~a_{j_1}^{[j_2]}\cdot a_{j_2}^{[j_3]}\cdot\, \cdots\, \cdot a_{j_{\alpha-2}}^{[j_{\alpha-1}]}\cdot a_{j_{\alpha-1}}^{[1]}=&\ \binom{n}{4}\cdot\sum a_5^{[j_1]}\cdot a_{j_1}^{[j_2]}\cdot a_{j_2}^{[j_3]}\cdot a_{j_3}^{[2]}\\
=&\ \binom{n}{4}\cdot\left(a_5^{[4]}\cdot a_4^{[3]}\cdot a_3^{[1]}\cdot a_2^{[\,1\,]}\right)\\
=&\ \binom{n}{4}\cdot\left(24a_2^4\right)~ {\bf \cdot}
\end{align*}
So, it follows that
\begin{align*}
f_5^{(n)}=&\ ~a_5\cdot\binom{n}{1}+\sum_{\alpha= 2}^{4}\left[ \binom{n}{\alpha}\cdot\sum a_5^{[j_1]}\cdot ~a_{j_1}^{[j_2]}\cdot a_{j_2}^{[j_3]}\cdot\, \cdots\, \cdot a_{j_{\alpha-2}}^{[j_{\alpha-1}]}\cdot a_{j_{\alpha-1}}^{[1]}\right]\\
&\\
= &\ a_5\cdot\binom{n}{1}+A_{2,5}+ A_{3,5}+ A_{4,5}~ {\bf \cdot}
\end{align*}
Thus,
\[f_5^{(n)}=~a_5\cdot\binom{n}{1}+\binom{n}{2}\cdot\left( ~5a_2^2+6a_2a_4+3a_3^2 ~\right)+\binom{n}{3}\cdot\left( ~10a_2^4+26a_2^2a_3~\right)+24~\binom{n}{4}\cdot a_2^4~~ {\bf \cdot}\]
\newpage
\bibliographystyle{amsplain}

\end {document}